%% file: mmcsurvey.tex
\documentclass[a4paper]{article}

\usepackage[utf8]{inputenc}
\usepackage[T1]{fontenc}
\usepackage{lmodern}
\usepackage{fullpage}
\usepackage{microtype}
\usepackage[usenames,dvipsnames,table]{xcolor}
\usepackage{amsmath,amsthm,amssymb}
\usepackage{wasysym} 
\usepackage{xspace}
\usepackage{subfig}
\usepackage{graphicx}
\usepackage{xifthen}
\usepackage{framed}
\usepackage{todonotes}
\usepackage{booktabs}
\usepackage{tabularx}
\usepackage{multirow}

\usepackage{tikz} 
\usetikzlibrary{matrix}
\usetikzlibrary{arrows,shapes,positioning}
\usetikzlibrary{decorations.pathreplacing}
\usetikzlibrary{calc}
 \tikzstyle{tredge}=[ultra thick, draw = nicered]

\usepackage[hidelinks,colorlinks,allcolors=blue]{hyperref}

\newtheorem{theorem}{Theorem} 
\newtheorem{proposition}{Proposition} 
\newtheorem{corollary}{Corollary}
\newtheorem{open}{Open Problem} 
\newtheorem{lemma}[theorem]{Lemma} 
\newcommand{\mc}{{\sc Matching Cut}}  
\newcommand{\mmc}{{\sc Maximum Matching Cut}}          
\newcommand{\minmc}{{\sc Minimum Matching Cut}}       
\newcommand{\maxcut}{{\sc Max Cut}}
\newcommand{\mincut}{{\sc Min Cut}}         
\newcommand{\NP}{{\sf NP}}
\newcommand{\ssi}{\subseteq_i}

\definecolor{nicered}{RGB}{204,0,0}
\definecolor{lightblue}{RGB}{153,204,255}                 
                           
\tikzstyle{bvertex}=[thin,circle,inner sep=0.cm, minimum size=1.7mm, fill=lightblue, draw=lightblue]
    \tikzstyle{rvertex}=[thin,circle,inner sep=0.cm, minimum size=1.7mm, fill=nicered,draw=nicered]
    \tikzstyle{vertex}=[thin,circle,inner sep=0.cm, minimum size=1.7mm, fill=black, draw=black]    
                 \tikzstyle{hedge}=[thin, draw = gray]
    \tikzstyle{edge} = [thick, draw=gray]
    \tikzstyle{cutedge}=[ultra thick, draw=nicered]
    \tikzstyle{tedge}=[ultra thick, draw=black]

\author{
  Van Bang Le\thanks{Universit\"at Rostock, Institut f\"ur Informatik, Germany -- \texttt{van-bang.le@uni-rostock.de}}
  \and Felicia Lucke\thanks{Department of Informatics, University of Fribourg, Fribourg, Switzerland -- \texttt{felicia.lucke@unifr.ch}}
  \and Daniël Paulusma\thanks{Department of Computer Science, Durham University, Durham, UK -- \texttt{daniel.paulusma@durham.ac.uk}}
  \and Bernard Ries\thanks{Department of Informatics, University of Fribourg, Fribourg, Switzerland -- \texttt{bernard.ries@unifr.ch}}
}

\title{Maximizing Matching Cuts}

\begin{document}
\maketitle

\begin{abstract}
A matching cut in a graph $G$ is an edge cut of $G$ that is also a matching.
This short survey gives an overview of old and new results and open problems for  {\sc Maximum Matching Cut}, which is to determine the size of a largest  matching cut in a graph. We also compare this problem with the related problems
{\sc Matching Cut}, {\sc Minimum Matching Cut}, and {\sc Perfect Matching Cut}, which are to determine if a graph has a matching cut; the size of a smallest matching cut in a graph; and if a graph has a matching cut that is a perfect matching, respectively. Moreover, we discuss a relationship between \mmc\ and \maxcut, 
which is to determine the size of a largest edge cut in a graph, as well as a relationship between {\sc Minimum Matching Cut} and \mincut, which is to determine the size of a smallest edge cut in a graph.
\end{abstract}

\section{Introduction}

Graph cut problems belong to a well-studied class of classical graph problems related to network connectivity, which is a central concept within theoretical computer science. 
More formally, in a graph $G=(V,E)$, a subset of edges $M\subseteq E$ is called an \emph{edge cut} if there exists a partition $(R,B)$ of $V$ into two non-empty subsets $R$ (say, of {\it red} vertices) and $B$ (say, of {\it blue} vertices) such that $M$ consists of exactly those edges whose end-vertices have different colours, so one of them belongs to $R$ and the other to $B$.

\begin{figure}[t]
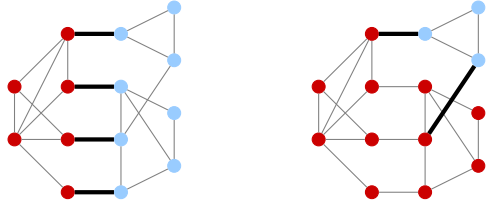

\centering
\include{defmmc}
\caption{\label{f-mcmmc}An example of a maximum matching cut (left) and a minimum matching cut (right) in a graph, where the edge cut is indicated in bold.}\label{f-maxmin}
\end{figure}

An edge cut $M$ in a graph~$G$ is a \emph{maximum}, respectively, \emph{minimum} edge cut of $G$ if $M$ has maximum, respectively, minimum size over all edge cuts in~$G$.
This leads to the two famous graph cut problems: \maxcut\ and \mincut, which are to determine for a given graph, the size of an edge cut with maximum, respectively, minimum number of edges.  
While \maxcut\  is a classical \NP-complete problem~\cite{GJS76}\footnote{The version of \maxcut\ with weights on edges is one of Karp's original 21 \NP-complete problems~\cite{Ka72}. In this survey we do not consider edge weights.}; 
(see also the surveys~\cite{Co09,PT95}), there exist several polynomial-time algorithms~\cite{FF56,Ga95} for \mincut\ (as \mincut\ can be modelled as a maximum flow problem).

In this short survey, we consider a notion of edge cuts which has received much attention lately. An edge cut $M$ of a graph is a \emph{matching cut} if $M$ is a {\it matching}, that is, no two edges of $M$ share an end-vertex, or equivalently, every red vertex has at most one blue neighbour, and vice versa. Not every graph has a matching cut (consider, for example, a triangle) and a graph may have multiple matching cuts (see Figure~\ref{f-maxmin}). Graphs with a matching cut were introduced by Graham~\cite{Gr70} in 1970 as {\it decomposable} graphs and were used to solve a problem on cube numbering~\cite{Gr70}. Matching cuts also have applications in ILFI networks~\cite{FP82}, graph drawing~\cite{PP01}, graph homomorphism problems~\cite{GPS12} and were used for determining conflict graphs for WDM networks~\cite{ACGH12}.
The decision problem {\sc Matching Cut} is to decide if a given graph has a matching cut. This problem was shown to be \NP-complete by Chv\'atal~\cite{Ch84}.

A matching cut $M$ is {\it maximum}, respectively, {\it minimum} if $M$ has maximum, respectively, minimum size over all matching cuts in~$G$ (if $G$ is decomposable). We refer to Figure~\ref{f-maxmin} for an example of a graph with a maximum and minimum matching cut of different size. We focus on the two corresponding optimization versions of \mc, which are the analogs of \maxcut\ and \mincut, respectively:

\medskip\noindent
\fbox{
\begin{minipage}{.965\textwidth}
\textsc{Maximum Matching Cut}\\[.7ex]
\begin{tabular}{l l}
{\em Instance:\/}& A graph $G$.\\
{\em Task:\/}& Determine the size of a maximum matching cut in $G$.
\end{tabular}
\end{minipage}
}

\medskip\noindent
\fbox{
\begin{minipage}{.965\textwidth}
\textsc{Minimum Matching Cut}\\[.7ex]
\begin{tabular}{l l}
{\em Instance:\/}& A graph $G$.\\
{\em Task:\/}& Determine the size of a minimum matching cut in $G$.
\end{tabular}
\end{minipage}
}

\bigskip
\noindent
Both \maxcut\ and \mincut\ are \NP-hard, as {\sc Matching Cut} is \NP-complete~\cite{Ch84}. The fact that \mmc\ is \NP-hard also follows from the \NP-completeness of another problem variant, known as {\sc Perfect Matching Cut}. A matching cut $M$ in a graph~$G$ is \emph{perfect} if $M$ is a {\it perfect} matching, that is, every vertex of $G$ is incident to an edge of $M$, or equivalently, every red vertex has exactly one blue neighbour, and vice versa. Any perfect matching cut of a graph is maximum, but the reverse might not be necessarily true. Heggernes and Telle~\cite{HT98} included perfect matching cuts in their $(\sigma,\rho)$-vertex partitioning framework
and proved that the corresponding decision problem, {\sc Perfect Matching Cut}, is \NP-complete.

\medskip
\noindent
{\bf Outline.}
The {\sc Maximum Matching Cut} problem was introduced in~\cite{LL23,LPR23b} and is our main focus. Due to its \NP-hardness, \mmc\ was studied for special graph classes. We survey these results in Section~\ref{s-mmc}. In the same section we also discuss 
 a result in~\cite{LPR23b} that shows how  {\sc Max Cut} can be reduced to {\sc Maximum Matching Cut}.
 Throughout Section~\ref{s-mmc} we compare known complexity results for special graph classes with corresponding results for {\sc Matching Cut} and {\sc Perfect Matching Cut}. In particular, we show how these results may differ from each other, which enables us to identify a number of open problems. In Section~\ref{s-minmc} we consider the {\sc Minimum Matching Cut} as a natural counterpart of \mmc, which was not studied before. 
Our aim in this section is to show, apart from some interesting open problems, that the complexities of \mmc\ and {\sc Minimum Matching Cut} may differ on special graph classes.
We do this by illustrating that for some graph classes, it is possible to reduce {\sc Minimum Matching Cut} to {\sc Min Cut}. 
 We conclude our survey with some final open problems in Section~\ref{s-con}.

\section{Preliminaries}\label{s-pre}

All graphs considered in this survey are undirected and have no multiple edges and self-loops, unless explicitly said otherwise.

Let $G=(V,E)$ be a graph.
A subset $S\subseteq V$ is a {\it clique} of $G$ if all vertices of $S$ are pairwise adjacent, whereas $S$ is an {\it independent set} if all vertices of $S$ are pairwise non-adjacent.
If $uv$ is an edge in $E$, then $u$ and $v$ are {\it neighbours}. For a vertex $u\in V$, the set $N(u)=\{v \in V\; |\; uv\in E\}$ denotes the {\it neighbourhood} of $u$.
The {\it degree} of $u$ in $G$ is the size $|N(u)|$ of the neighbourhood of $u$. 
If every vertex of $G$ has degree~$r$ for some integer $r\geq 0$, we say that $G$ is {\it $r$-regular}. If we just say that $G$ is {\it regular}, we mean that there exists an integer $r$ such that $G$ is $r$-regular.
We say that $G$ is \emph{subcubic} if every vertex of $G$ has degree at most~$3$ and that $G$ is \emph{cubic} if every vertex of $G$ has degree exactly~$3$.
The {\it line graph} of $G$ is the graph that has vertex set~$E(G)$, such that there is an edge between two vertices $e_1$ and $e_2$ if and only if  $e_1$ and~$e_2$ share an end-vertex in $G$.

A graph $G=(V,E)$ is {\it bipartite} if $V$ can be partitioned into two independent sets $A$ and $B$, which are called the {\it partition classes} of $G$. A bipartite graph $G$ is {\it $(k,\ell)$-regular} if it has partition classes $A$ and $B$, such that every vertex in $A$ has degree~$k$ in $G$, and every vertex in $B$ has degree $\ell$ in~$G$.

The {\it length} of a graph that is either a path or a cycle is its number of edges. 
A cycle is said to be {\it odd} if it has odd length. The \emph{girth} of a graph~$G$ is the length of a shortest cycle of $G$. If $G$ is a {\it forest} (that is, a graph with no cycles), then $G$ has girth $\infty$.

Let $G=(V,E)$ be a graph. The {\it distance} between two vertices $u$ and $v$ in~$G$ is the {\it length} of a shortest path between $u$ and $v$ in $G$. The {\it eccentricity} of $u$ is defined as the maximum distance between $u$ and any other vertex of $G$.  This gives us the {\it diameter} of $G$, which is the maximum eccentricity over all vertices of~$G$, and the {\it radius} of $G$, which is the minimum eccentricity over all vertices of~$G$.  We note that the radius of $G$ is at most its diameter, whereas the diameter of $G$ is at most twice its radius.

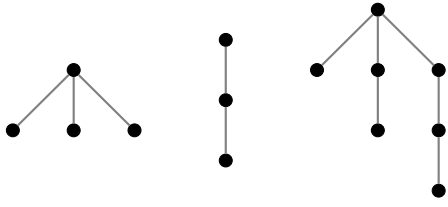
\begin{figure}[t]
\centering
\input{s-graphs}
\caption{The graph $K_{1,3}+P_3+S_{1,23}$, which is an example of a graph that belongs to $\mathcal{S}$.}\label{f-s-graphs}
\end{figure}

The graphs $P_r$, $C_s$, $K_t$ denote the path, cycle and complete graph on $r$, $s$, and $t$ vertices, respectively. The {\it diamond} is the graph obtained from $K_4$ after removing an edge.
The graph $K_{1,\ell}$ denotes the {\it star} on $\ell+1$ vertices, which is the (bipartite) graph on vertices $u,v_1,\ldots,v_\ell$ with edges $uv_i$ for every $i\in \{1,\ldots,\ell\}$.  The graph $K_{1,3}$ is commonly known as the {\it claw}.
For $1\leq h\leq i\leq j$, the graph $S_{h,i,j}$ is the tree with one vertex $u$ of degree~$3$,
whose (three) leaves are at distance~$h$,~$i$ and~$j$ from $u$. We note that $S_{1,1,1}=K_{1,3}$. We say that a graph $S_{h,i,j}$ is a {\it subdivided claw}.
The {\it disjoint union} $G_1+G_2$ of two graphs $G_1$ and $G_2$ with $V(G_1)\cap V(G_2)=\emptyset$ is the graph $(V(G_1)\cup V(G_2),E(G_1)\cup E(G_2))$. We denote the disjoint union of $s$ copies of the same graph $G$ by $sG$.
We define the set ${\cal S}$ as all graphs that are the disjoint union of one or more graphs, each of which is either a path or a subdivided claw, see Figure~\ref{f-s-graphs} for an example.

For a set $S\subseteq V(G)$, the graph $G[S]$ is the subgraph of a graph $G$ {\it induced} by $S$, which is the graph obtained from $G$ after deleting every vertex that does not belong to $S$.
For an integer $\ell$, a graph $G$ on more than $\ell$ vertices is said to be {\it $\ell$-connected} if $G[V\setminus S]$ is connected for every set $S$ on at most $\ell-1$ vertices. 
An edge $e$ of a connected graph $G$ is a {\it bridge} if the graph $G-e$, obtained from $G$ after deleting~$e$, is disconnected; note that an edge is a bridge if and only if $M=\{e\}$ is a matching cut.
A graph with no bridges is said to be {\it bridgeless}.

Let $G$ and $H$ be two graphs. We say that $G$ contains $H$ as an {\it induced subgraph} if $G$ can be modified to $H$ by a sequence of vertex deletions; if not, then $G$ is said to be {\it $H$-free}. 
We use the notation $H\ssi G$ to indicate that $H$ is an induced subgraph of $G$.
We say that $G$ contains $H$ as a {\it subgraph} if $G$ can be modified to $H$ by a sequence of vertex deletions and edge deletions; if not, then $G$ is {\it $H$-subgraph-free}. 
Moreover, $G$ contains $H$ as a  {\it spanning subgraph} if $G$ can be modified to $H$ by a sequence of only edge deletions (so $V(G)=V(H)$).

In order to define some more graph containment relations, we first need to define some more graph operations. Let $G=(V,E)$ be a graph.
The {\it contraction} of an edge $e=uv$ in~$G$ replaces $u$ and $v$ by a new vertex $w$ that is made adjacent to every vertex of $(N(u)\cup N(v))\setminus \{u,v\}$  (without creating multiple edges, 
unless we explicitly say otherwise). 
Suppose that one of $u,v$, say $v$, had degree~$2$ in $G$, and moreover that the two neighbours of $v$ in $G$ are not adjacent. In that case, we also say that by contracting $uv$, we {\it dissolved} $v$, and in this specific case we also call the edge contraction a \emph{vertex dissolution}, namely of vertex $v$.
For a subset $S\subseteq V$, we say that we {\it contract $G[S]$ to a single vertex} in $G$ if we contract every edge of a spanning tree of $G[S]$. 

Let $G$ and $H$ be two graphs.
We say that $G$ contains $H$ as a {\it topological minor} (or as a {\it subdivision}) if $G$ can be modified to $H$ by a sequence of vertex deletions, vertex dissolutions and edge deletions; if not, then $G$ is {\it $H$-topological-minor-free}. Likewise, $G$ contains $H$ as a {\it minor} if $G$ can be modified to $H$ by a sequence of vertex deletions, edge deletions and edge contractions; if not, then $G$ is {\it $H$-minor-free}. 

Finally, let $G$ be a graph and ${\cal H}$ be a set of graphs. We say that $G$ is {\it ${\cal H}$-free} if $G$ is $H$-free for every $H\in {\cal H}$. We define the notions of being {\it ${\cal H}$-subgraph-free, ${\cal H}$-topological-minor-free and ${\cal H}$-minor-free} analogously. If ${\cal H}=\{H_1,\ldots,H_p\}$ for some integer $p\geq 1$, then we may also write that $G$ is $(H_1,\ldots,H_p)$-free. We say that $G$ is  {\it quadrangulated} if $G$ is ${\cal C}_{\geq 5}$-free, where ${\cal C}_{\geq 5}= \{C_5,C_6,\ldots\}$.

\section{Maximum Matching Cut}\label{s-mmc}

 \begin{figure}[t]
  \vspace*{5mm}
 \begin{center}
 \begin{tikzpicture}[scale=1]
\node[] (g) at (-1,0.5){$G$};

\node[vertex] (v1) at (0,0){};
\node[vertex] (v2) at (1,0){};
\node[vertex] (v3) at (0,1){};
\node[vertex] (v4) at (1,1){};

\draw[tredge] (v1)--(v2);
\draw[tredge] (v3)--(v2);
\draw[edge] (v4)--(v2);
\draw[tredge] (v3)--(v4);

\begin{scope}[shift = {(3,-0.5)}, scale = 0.5]
\node[] (g) at (6,2){$G'$};

\node[vertex] (v11) at (1,0){};
\node[vertex] (v12) at (1,1){};
\node[vertex] (v13) at (0,1){};

\node[vertex] (v21) at (3,0){};
\node[vertex] (v22) at (3,1){};
\node[vertex] (v23) at (4,1){};

\node[vertex] (v31) at (0,3){};
\node[vertex] (v32) at (1,3){};
\node[vertex] (v33) at (1,4){};

\node[vertex] (v41) at (4,3){};
\node[vertex] (v42) at (3,3){};
\node[vertex] (v43) at (3,4){};

\draw[edge] (v11)--(v12);
\draw[edge] (v12)--(v13);
\draw[edge] (v13)--(v11);

\draw[edge] (v21)--(v22);
\draw[edge] (v22)--(v23);
\draw[edge] (v23)--(v21);

\draw[edge] (v31)--(v32);
\draw[edge] (v32)--(v33);
\draw[edge] (v33)--(v31);

\draw[edge] (v41)--(v42);
\draw[edge] (v42)--(v43);
\draw[edge] (v43)--(v41);

\draw[tredge] (v11)--(v21);
\draw[tredge] (v22)--(v32);
\draw[edge] (v23)--(v41);
\draw[tredge] (v43)--(v33);

\end{scope}
\end{tikzpicture}
 \caption{The example from~\cite{LPR23b} that illustrates the polynomial-time reduction from \maxcut\ to \mmc. Left: a connected graph $G$ with $\Delta=3$, where the thick red edges form a maximum edge cut. Right: the corresponding graph~$G'$, where the thick red edges form a maximum matching cut.}\label{fig-clawfree}
 \end{center}
 \end{figure}
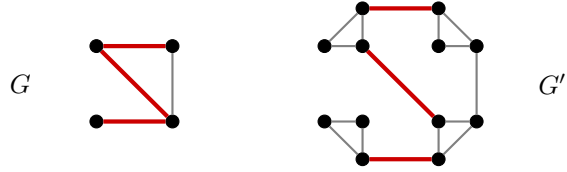

We first discuss, in Section~\ref{s-reduc}, a strong relationship between \mmc\ and \maxcut. Afterwards we focus, in Section~\ref{s-para}, on results for classes of graphs where some distance metric or connectivity parameter is bounded. Finally, in Section~\ref{s-contain}, we consider \mmc\ for graph classes defined by some containment relation. 

\subsection{Reducing Max Cut to Maximum Matching Cut}\label{s-reduc}

We first observe that the complexities of  \mmc\ and \maxcut\ may differ on special graph classes. For example, {\sc Max Cut} is polynomial-time solvable for planar 
graphs~\cite{Ha75} and trivial for bipartite graphs,
whereas \mmc\ (even \mc) remains \NP-hard when restricted to planar graphs of 
girth~$5$~\cite{Bo09} and bipartite graphs of maximum degree~$4$~\cite{LR03}.
However, there exists a simple polynomial-time reduction from \maxcut\ to \mmc\ provided in~\cite[Theorem 20]{LPR23b}, which has some interesting consequences.  

Given a connected graph $G$ with maximum degree $\Delta\ge 3$, let $G'$ be the graph obtained from $G$ as follows (see also Figure \ref{fig-clawfree}):
\begin{itemize}
\item [1.] replace each vertex $v$ of $G$ by a clique $C_v$ on $\Delta$ vertices;
\item [2.] for each edge $uv$ of $G$, add an edge between a vertex in $C_u$ and a vertex in $C_v$ such that for every vertex $w$ of $G$, every vertex in $C_w$ has at most one neighbour outside $C_w$. 
\end{itemize}

\noindent
We can now make the following observation, shown in~\cite{LPR23b} as part of the proof of Theorem~20, except that it was not observed in~\cite{LPR23b} that $G'$ has
no induced odd cycle of length at least~$5$.

\begin{proposition}\label{p-gaccent}
A connected graph $G$ 
with maximum degree~$\Delta\geq 3$ 
has an edge cut of size 
$k$ if and only if the corresponding graph $G'$ has a matching cut of size 
$k$. Moreover, the following holds:
\begin{itemize}
\item $G'$ has maximum degree $\Delta$; 
\item $G'$ is regular if $G$ is regular (of degree~$\Delta$);
\item $G'$ is $(\text{claw},\text{diamond}\/)$-free and has no induced cycle of odd length at least~$5$, that is, $G'$ is a line graph of a bipartite graph. 
\end{itemize}
\end{proposition}

\noindent
It is well-known that \maxcut\ can be approximated within a ratio of $0.878567$~\cite{GW95}. However, assuming the Unique Games Conjecture, it is \NP-hard to approximate \maxcut\ within a ratio better than $0.878567$~\cite{KKMO07}. For cubic graphs, \maxcut\ admits an approximation within a better ratio of~$0.9326$~\cite{HalperinLZ02} but is still {\sf APX}-hard~\cite{AK00}. 
Hence, Proposition~\ref{p-gaccent} has the following two consequences, the first of which is stronger than the \NP-hardness result in~\cite[Theorem~20]{LPR23b}.

\begin{corollary}\label{c-two}
{\sc Maximum Matching Cut} is {\sf APX}-hard even for cubic line graphs of bipartite graphs, and assuming the Unique Games Conjecture, \NP-hard to approximate within a ratio better than $0.878567$.
\end{corollary}

\noindent
An {approximation with ratio $c$} for \mmc\ is a polynomial-time algorithm that outputs either {\tt no} if the input graph $G$ has no matching cut, or else a matching cut of size at least $c\cdot{\sf opt}$, where {\sf opt} is the maximum size of a matching cut in $G$. We pose the following open problem.

\begin{open}\label{o-c}
Is there an approximation for \mmc\ within some constant ratio when restricted to (sub)cubic graphs?
\end{open}

\noindent
Note that it follows from the {\sf NP}-completeness of \mc\ that in the class of all graphs, there exists no approximation for \mmc\ within any ratio. 

\subsection{Bounding the Degree, Diameter, Radius or Girth}\label{s-para}

In this section, we will present dichotomies for \mmc\ restricted to graphs of bounded maximum degree, diameter, radius or girth. We will also illustrate how these dichotomies differ from corresponding complexity results for {\sc Matching Cut} and {\sc Perfect Matching Cut}.

We first consider the maximum degree of a graph.
It is known that {\sc Matching Cut} is polynomial-time solvable for 
 subcubic graphs~\cite{Ch84} but \NP-complete even for
$(3,4)$-regular bipartite graphs~\cite{LR03}.
Bonnet, Chakraborty and Duron~\cite{BCD23} proved that {\sc Perfect Matching Cut} is \NP-complete for $3$-connected cubic planar bipartite graphs. Combining these results with Corollary~\ref{c-two} yields the following result.

 \begin{theorem}\label{t-degree}
 \mmc\ is \NP-hard for  $(3,4)$-regular bipartite graphs; $3$-connected cubic planar bipartite graphs; and cubic line graphs of bipartite graphs.
 \end{theorem}
 
\noindent
Theorem~\ref{t-degree}, combined with the straightforward observation that \mmc\ is polynomial-time solvable for graphs of maximum degree at most~$2$, leads to the following dichotomy, which is different from the above dichotomy for {\sc Matching Cut}.
 
 \begin{corollary}\label{c-dichodegree}
For an integer~$\Delta$,  \mmc\ on (bipartite) graphs of maximum degree at most $\Delta$ is  polynomial-time solvable if $\Delta \leq 2$, and \NP-hard if $\Delta \geq 3$.
\end{corollary}

 \noindent
 We now consider graphs of bounded diameter and graphs of bounded radius.
 The {\sc Matching Cut} problem is polynomial-time solvable for graphs of radius (and thus diameter) at most~$2$~\cite{LPR22}, but \NP-complete for graphs of diameter
 (and thus radius)
  $3$~\cite{LL19}.
 The {\sc Perfect Matching Cut} problem is polynomial-time solvable for graphs of radius (and thus diameter) at most~$2$~\cite{LPR23a}, but \NP-complete for graphs of diameter~$4$
 and 
 radius~$3$~\cite{LL23}.\footnote{The latter result follows from the \NP-hardness gadget in~\cite{LL23} for {\sc (Perfect) Matching Cut} for $(3P_6,2P_7,P_{14})$-free graphs. The complexity status of {\sc Perfect Matching Cut} for graphs of diameter~$3$ is not known.} 
 However, for \mmc, the following result is known.
 
 \begin{theorem}[\cite{LPR23b}]\label{t-quadra}
\mmc\  is polynomial-time solvable for graphs of diameter at most~$2$ but \NP-hard for $2P_3$-free quadrangulated graphs of diameter~$3$ and radius~$2$.
\end{theorem}

\noindent
Theorem~\ref{t-quadra} yields the following two dichotomies, which show that from the above results, only the results for {\sc Matching Cut} on graphs of bounded diameter can be generalized to {\sc Maximum Matching Cut}.

\begin{corollary}[\cite{LPR23b}]\label{c-dichodiam}
For an integer~$d$, \mmc\ on graphs of diameter $d$ is  polynomial-time solvable if $d \leq 2$, and \NP-hard if $d \geq 3$.
\end{corollary}

\begin{corollary}[\cite{LPR23b}]\label{c-dichorad}
For an integer~$r$, \mmc\ on graphs of radius $r$~is polynomial-time solvable if $r \leq 1$, and \NP-hard if $r \geq 2$.
\end{corollary}

\noindent
It is known that {\sc Matching Cut} is polynomial-time solvable for bipartite graphs of diameter at most~$3$ and \NP-complete for bipartite graphs of diameter~$4$~\cite{LL19}. 
The latter result implies that \mmc\ is \NP-complete for bipartite graphs of diameter~$4$. We pose the following open problem: 

\begin{open}\label{o-b3}
Determine the complexity of \mmc\ for bipartite graphs of diameter at most~$3$.
\end{open}

\noindent
To solve Open Problem~\ref{o-b3}, we cannot use the \NP-hardness gadgets from~\cite{BCD23,LR03} as these gadgets are used to prove \NP-completeness for bipartite graph classes of bounded degree (see Theorem~\ref{t-degree}) and thus naturally have unbounded diameter (in particular, gadgets of bounded degree and bounded diameter would have bounded size).

We now consider the girth. Recently, Feghali et al.~\cite{FLPR23} proved that for every integer $g \geq 3$, {\sc Matching Cut} is \NP-complete for bipartite graphs of girth at least~$g$ and maximum degree at most~$60$. This immediately leads to the following result.

\begin{theorem}\label{t-dichogirth}
For every integer~$g\geq 3$, \mmc\ is \NP-hard for graphs of girth at least~$g$ and maximum degree at most~$60$.
\end{theorem}

\noindent
For \mmc\ it might be possible to find an alternative hardness gadget, and we pose the following open problem:

\begin{open}
Is it possible to improve the bound on the maximum degree in Theorem~\ref{t-dichogirth}?
\end{open}

\subsection{Graph Containment Relations}\label{s-contain}

In this section we consider graph classes defined by some containment relation. We consider forbidden induced subgraphs, minors, topological minors and subgraphs.

We start with the induced subgraph relation. We first consider $H$-free graphs for some graph $H$. The complexity classifications for {\sc Matching Cut} and {\sc Perfect Matching Cut} are still open; so far all known results for specific graphs $H$ suggest that there is no graph~$H$, such that these two problems differ in complexity when restricted to $H$-free graphs. However, for both problems, we must still solve a large number of cases where $H$ is a disjoint union of paths and subdivided claws 
(see, for example,~\cite{FLPR23}). 

In contrast to the above, for \mmc, we have a complete classification of its complexity for $H$-free graphs, as shown in~\cite{LPR23b}. 
Namely, to obtain \NP-hardness, we can apply Theorem~\ref{t-dichogirth} if $H$ contains a cycle; Theorem~\ref{t-degree} if $H$ contains an induced claw\footnote{The class of line graphs is readily seen to be contained in the class of claw-free graphs.}; and Theorem~\ref{t-quadra} if $H$ contains an induced~$2P_3$. In all remaining cases, $H$ is an induced subgraph of $sP_2+P_6$ for some $s\geq 0$, for which the problem is polynomial-time solvable~\cite{LPR23b}. 

\begin{theorem}[\cite{LPR23b}]\label{t-dichoH}
For a graph~$H$, \mmc\ on $H$-free graphs is 
polynomial-time solvable if $H\ssi sP_2+P_6$ for some $s\geq 0$, and \NP-hard otherwise.
\end{theorem}

\noindent
We note that there exist graphs~$H$ for which the complexity of \mmc\ is different from the complexity of {\sc (Perfect) Matching Cut} (subject to ${\sf P}\neq \NP$); 
in particular, {\sc Matching Cut} and {\sc Perfect Matching Cut} are polynomial-time solvable for $(sP_3+P_6)$-free graphs~\cite{LPR22} and $(sP_4+P_6)$-free graphs~\cite{LPR23a}, respectively, 
for any integer $s\geq 1$, while \mmc\ is \NP-hard even for $2P_3$-free graphs by Theorem~\ref{t-dichoH}.

We pose the following extension of Theorem~\ref{t-dichoH} as an open problem.

\begin{open}
For every finite set of graphs~${\cal H}$, determine the complexity of \mmc\ for ${\cal H}$-free graphs.
\end{open}

\noindent
We now consider ${\mathcal H}$-minor-free graphs and ${\mathcal H}$-topological-minor-free graphs. 
As a well-known consequence of a classic result of~Robertson and Seymour~\cite{RS86}\footnote{For an explicit explanation of this consequence of~\cite{RS86}, see, for example~\cite{JMOPPSV,KLMT11}.}, any graph problem~$\Pi$ that is \NP-hard on subcubic planar graphs but polynomial-time solvable for graphs of bounded treewidth can be fully classified on ${\cal  H}$-minor-free graphs and ${\cal H}$-topological minor-free graphs, even for infinite sets~${\cal H}$. Namely, $\Pi$ on ${\cal H}$-minor-free graphs is polynomial-time solvable if ${\cal H}$ contains a planar graph and \NP-hard otherwise, while $\Pi$ on ${\cal H}$-topological-minor-free graphs is polynomial-time solvable if ${\cal H}$ contains a subcubic planar graph and \NP-hard otherwise. 

It follows from the framework of Arnborg, Lagergren and Seese~\cite{ALS91} that \mmc\ is polynomial-time solvable for graphs of bounded treewidth. Hence, combining this observation with the above results of~\cite{RS86} and Theorem~\ref{t-degree} yields the following two dichotomies.

\begin{theorem}\label{t-planar} 
For any set of graphs~${\cal H}$, \mmc\ on ${\cal H}$-(topological-)minor-free graphs is polynomial-time solvable if ${\cal H}$ contains a (subcubic) planar graph and is \NP-hard otherwise.
\end{theorem}

\noindent
Finally, we consider ${\mathcal H}$-subgraph-free graphs. 
Recall that ${\cal S}$ denotes the class of graphs that are disjoint unions of paths and subdivided claws and that \mmc\ is polynomial-time solvable for graphs of bounded treewidth. The latter implies that \mmc\ is polynomial-time solvable for ${\cal H}$-subgraph-free graphs if ${\cal H}$ contains a graph from ${\cal S}$~\cite{RS84}.
In~\cite{FLPR23}, results from~\cite{LT22} were combined to prove that
for any finite set of graphs ${\cal H}$, {\sc Perfect Matching Cut} is \NP-complete for ${\cal H}$-subgraph-free graphs if ${\cal H}$ contains no graph from ${\cal S}$.
Hence, we find the following dichotomy for \mmc\ (which is the same dichotomy as for {\sc Perfect Matching Cut}~\cite{FLPR23}).

\begin{theorem}\label{t-dichosubgraph}
For any finite set of graphs ${\cal H}$, {\sc Maximum Matching Cut} on ${\cal H}$-subgraph-free graphs is polynomial-time solvable if ${\cal H}$ contains a graph from ${\cal S}$ and \NP-hard otherwise.
\end{theorem}

\noindent
In contrast to Theorem~\ref{t-planar}, we note that Theorem~\ref{t-dichosubgraph} only holds for {\it finite} sets of graphs~${\cal H}$, and we refer to~\cite{JMOPPSV} for examples that show that this condition cannot be avoided. The following problem looks challenging.

\begin{open}
Classify the complexity of \mmc\ for ${\cal H}$-subgraph-free graphs if ${\cal H}$ is infinite.
\end{open}

\section{Minimum Matching Cut}\label{s-minmc}

In this section we focus on {\sc Minimum Matching Cut}. We first show that there exists a graph class, for which the complexities of {\sc Maximum Matching Cut} and {\sc Minimum Matching Cut} are different. Namely, we consider the class of claw-free graphs. By Theorem~\ref{t-dichoH}, {\sc Maximum Matching Cut} is \NP-complete for claw-free graphs. However, we show that the polynomial-time algorithm of Bonsma~\cite{Bo09} for {\sc Matching Cut} on claw-free graphs
can be extended to work for {\sc Minimum Matching Cut} via a reduction to {\sc Min Cut}.

A (connected) component of a graph $G$ is {\it non-trivial} if it has at least two vertices. Let $F$ be the set of edges of $G$ not in any triangle. Let $G[F]$ be the graph obtained from $G$ by deleting all edges not in $F$ and all vertices not incident to an edge of $F$. Let $G-F$ be the graph obtained from $G$ by removing the edges of $F$. 
We need the following lemma, which was proven by Bonsma; Lemma~\ref{l-bclaw}:2 is~\cite[Theorem~7]{Bo09} once 
we assume that $G$ is a bridgeless graph that is not a cycle, 
whereas Lemma~\ref{l-bclaw}:1 can be found within the proof of~\cite[Theorem~7]{Bo09}.

\begin{figure}[t]
\centering
\input{f-bonsma-clawfree}
\caption{Left: a claw-free graph $G_1$ that satisfies condition~(i) of Lemma~\ref{l-bclaw} with a minimum matching cut in bold. Middle: a  claw-free graph $G_2$ that satisfies condition~(ii) but not~(i) with a minimum matching cut in bold. Right: the (multi)graph $G_2^*$ with corresponding minimum edge cut in bold.}\label{f-bonsma-clawfree} 
\end{figure}
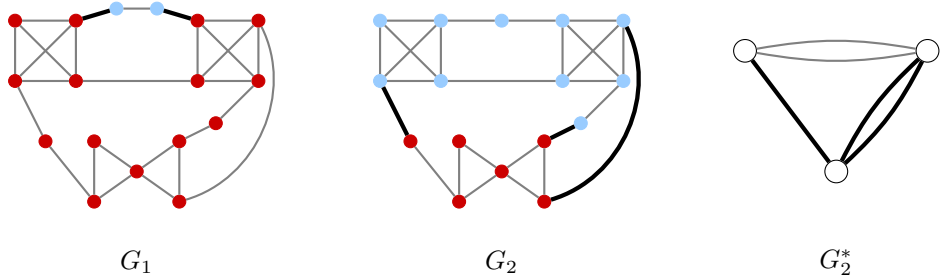

\begin{lemma}[\cite{Bo09}]\label{l-bclaw}
For a connected 
bridgeless 
claw-free graph $G$ that is not a cycle, 
 the following holds:
\begin{itemize}
\item [1.] Each component of $G[F]$ is a path (of length at least~$1$) in which every inner vertex has degree~$2$~in~$G$.
\item [2.] The graph $G$ has a matching cut if and only if
\begin{itemize}
\item [(i)] $G[F]$ contains a path component of length at least~$3$; or
\item [(ii)] $G-F$ contains at least two non-trivial components.
\end{itemize}
\end{itemize}
\end{lemma}

\noindent
Now, let $G$ be a connected claw-free graph. We may assume without loss of generality that $G$ is not a cycle and has no bridge; the latter implies that every matching cut of $G$ has size at least~$2$.
Conditions (i) and (ii) of Lemma~\ref{l-bclaw}:2 can be checked in polynomial time. If~(i) holds, then $G$ has a minimum matching cut of size~$2$ due to Lemma~\ref{l-bclaw}:1; see Figure~\ref{f-bonsma-clawfree}.
Now suppose (i) does not hold, but condition~(ii) holds. 
As~(i) does not hold, every component of $G[F]$ is a path of length~$1$ or~$2$. From the definition of $F$ it follows that the end-vertices of every path component of $G[F]$ of length~$2$ are not adjacent.
Moreover, as $G$ is bridgeless, the end-vertices of every path component of $G[F]$ belong to either the same non-trivial component of $G-F$ or to different non-trivial components of $G-F$.
The above means we can construct a new graph $G^*$ in polynomial time from~$G$ by

\begin{itemize}
\item [1.] dissolving every middle vertex of every path of length~$2$ in $G[F]$ whose end-vertices belong to different (non-trivial) connected components of $G-F$; 
\item [2.] deleting every middle vertex of every path of length~$2$ in $G[F]$ whose end-vertices belong to the same (non-trivial) connected component of $G-F$, and
\item [3.] finally, contracting every (non-trivial) component of $G-F$ to a single vertex without removing multiple edges.
\end{itemize}

\noindent
Note that $G^*$ may have multiple edges, as illustrated in Figure~\ref{f-bonsma-clawfree}. We now observe that any matching cut of size $k$ in $G$ corresponds to an edge cut of size $k$ in $G^*$, and vice versa. Indeed, first consider a matching cut $M$ of size $k$ in $G$. As no edge of a triangle can belong to a matching cut, $M$ can only contain edges of $G[F]$. Furthermore, since every component of $G[F]$ is a path of length~$1$ or~$2$, $M$ cannot have more than one edge of the same component of $G[F]$. It follows from the construction of $G^*$ that every path in $G[F]$ corresponds to an edge in $G^*$. Thus, for every edge of $M$ in $G$, we can choose a corresponding edge in $G^*$ to obtain an edge cut of size $k$ in $G^*$. Now conversely, assume $M^*$ is an edge cut of size $k$ in $G^*$. By construction, every edge in $G^*$ corresponds to a path in $G[F]$. Moreover, every subset of $E(G)$ that consists of at most one edge from each component of $G[F]$ is a matching in~$G$. Hence, by choosing for every edge in $M^*$ in $G^*$, an edge from the corresponding path component in $G[F]$, we obtain a matching cut of size $k$ in $G$. So, 
in order to solve \minmc{} on $G$, it suffices to solve {\sc Min Cut} on~$G^*$.  As the latter takes polynomial time~\cite{FF56}, we thus obtained the following result.

\begin{theorem}\label{t-mclaw}
\minmc{} is polynomial-time solvable for claw-free graphs.
\end{theorem}

\noindent
As mentioned, hardness results for {\sc Matching Cut} on special graph classes hold for {\sc Minimum Matching Cut}. However, it is still open which of the polynomial-time results carry over. In particular we ask:

\begin{open}\label{o-2p3}
Determine the complexity of {\sc Minimum Matching Cut} for $2P_3$-free graphs.
\end{open}

\noindent
We believe Open Problem~\ref{o-2p3} is interesting. 
On one hand, the \NP-hardness reduction from~\cite{LPR23b} for {\sc Maximum Matching Cut} on $2P_3$-free graphs fails for {\sc Minimum Matching Cut}.
On the other hand, the polynomial-time result from~\cite{LPR22} for {\sc Matching Cut} on $2P_3$-free graphs is based on a reduction to $2$-SAT and cannot be used either.

As in the case of \mmc, it follows from the {\sf NP}-completeness of \mc\ that there exists no polynomial-time approximation for \minmc.

\section{Summary}\label{s-con}

In this survey, we reviewed algorithmic and hardness results for \mmc, the maximization version of the classical decision problem \mc. We also pointed out that the complexities of \mmc\ and \minmc, the minimization version of \mc, may differ on special graph classes, and we proposed some relevant 
open problems for further research. 

To conclude, let us remark that both \mmc\ and \minmc\ are {\sf NP}-hard even in the following promise setting: the input graphs are given with the promise that every minimum matching cut is a maximum matching cut (in fact, every matching cut is a perfect matching cut, see~\cite{LL23}). 
In particular, recognizing graphs in which every maximal (minimal) matching cut is maximum (minimum) is {\sf NP}-complete. Thus, it would be interesting to characterize and recognize special graphs with this property. More precisely, we propose the following problem.

\begin{open}\label{o-min-max} 
Determine polynomially recognizable classes of graphs in which every maximal (minimal) matching cut is maximum (minimum).
\end{open}

\noindent
Examples of (non-trivial) graph classes in Problem~\ref{o-min-max} include 
the classes of $d$-dimensional hypercubes, $d\ge 2$. 
To make the problem more attractive, we point out that Lesk, Plummer and Pullblank~\cite{LPP84} proved that graphs in which every maximal matching is maximum can be recognized in polynomial time, while Chv\'atal and Slater~\cite{CS93} proved that it is {\sf coNP}-complete to recognize graphs in which every maximal independent set is maximum. 

\medskip
\noindent
{\it Acknowledgments.} The third author thanks Erik Jan van Leeuwen for a fruitful discussion on graph containment relations.

\bibliographystyle{plain}
\bibliography{../ref}

\end{document}

%% file: defmmc.tex
\begin{tikzpicture}

 \def\k{0.7}

	\node[rvertex] (a1) at (0,\k/2){};
	\node[rvertex] (a2) at (0,-\k/2){};
	\node[rvertex] (a3) at (\k,1.5*\k){};
	\node[rvertex] (a4) at (\k,\k/2){};
	\node[rvertex] (a5) at (\k,-3*\k/2){};
	\node[rvertex] (a5a) at (\k,-\k/2){};
	\node[bvertex] (a6) at (2*\k,1.5*\k){};
	\node[bvertex] (a7) at (2*\k,0.5*\k){};
	\node[bvertex] (a8a) at (2*\k,-0.5*\k){};
	\node[bvertex] (a8) at (2*\k,-1.5*\k){};
	\node[bvertex] (a9) at (3*\k,2*\k){};
	\node[bvertex] (a10) at (3*\k,1*\k){};
	\node[bvertex] (a11) at (3*\k,0*\k){};
	\node[bvertex] (a12) at (3*\k,-1*\k){};
	
	\draw[hedge](a1)--(a2);
	\draw[hedge](a1)--(a3);
	\draw[hedge](a1)--(a5a);
	\draw[hedge](a2)--(a3);
	\draw[hedge](a2)--(a4);
	\draw[hedge](a2)--(a5);
	\draw[hedge](a2)--(a5a);
	\draw[hedge](a3)--(a4);
	\draw[tedge](a3)--(a6);
	\draw[tedge](a4)--(a7);
	\draw[tedge](a5)--(a8);
	\draw[tedge](a5a)--(a8a);
	\draw[hedge](a6)--(a9);
	\draw[hedge](a6)--(a10);
	\draw[hedge](a7)--(a8a);
	\draw[hedge](a7)--(a11);
	\draw[hedge](a7)--(a12);
	\draw[hedge](a8)--(a12);
	\draw[hedge](a8a)--(a8);
	\draw[hedge](a9)--(a10);
	\draw[hedge](a11)--(a12);
         \draw[hedge](a8a)--(a10);
	
	\begin{scope}[shift= {(4,0)}]
	\node[rvertex] (a1) at (0,\k/2){};
	\node[rvertex] (a2) at (0,-\k/2){};
	\node[rvertex] (a3) at (\k,1.5*\k){};
	\node[rvertex] (a4) at (\k,\k/2){};
	\node[rvertex] (a5) at (\k,-3*\k/2){};
	\node[rvertex] (a5a) at (\k,-\k/2){};
	\node[bvertex] (a6) at (2*\k,1.5*\k){};
	\node[rvertex] (a7) at (2*\k,0.5*\k){};
	\node[rvertex] (a8a) at (2*\k,-0.5*\k){};
	\node[rvertex] (a8) at (2*\k,-1.5*\k){};
	\node[bvertex] (a9) at (3*\k,2*\k){};
	\node[bvertex] (a10) at (3*\k,1*\k){};
	\node[rvertex] (a11) at (3*\k,0*\k){};
	\node[rvertex] (a12) at (3*\k,-1*\k){};
	
	\draw[hedge](a1)--(a2);
	\draw[hedge](a1)--(a3);
	\draw[hedge](a1)--(a5a);
	\draw[hedge](a2)--(a3);
	\draw[hedge](a2)--(a4);
	\draw[hedge](a2)--(a5);
	\draw[hedge](a2)--(a5a);
	\draw[hedge](a3)--(a4);
	\draw[tedge](a3)--(a6);
	\draw[hedge](a4)--(a7);
	\draw[hedge](a5)--(a8);
	\draw[hedge](a5a)--(a8a);
	\draw[hedge](a6)--(a9);
	\draw[hedge](a6)--(a10);
	\draw[hedge](a7)--(a8a);
	\draw[hedge](a7)--(a11);
	\draw[hedge](a7)--(a12);
	\draw[hedge](a8)--(a12);
	\draw[hedge](a8a)--(a8);
	\draw[hedge](a9)--(a10);
	\draw[hedge](a11)--(a12);
         \draw[tedge](a8a)--(a10);
	
	\end{scope}

\end{tikzpicture}

%% file: s-graphs.tex
\begin{tikzpicture}
\begin{scope}[scale = 0.8]
\begin{scope}
\node[vertex](v1) at (0,2){};
\node[vertex](v2) at (0,1){};
\node[vertex](v3) at (0,0){};

\draw[edge](v1)--(v2);
\draw[edge](v2)--(v3);
\end{scope}

\begin{scope}[shift= {(2.5,0.5)}]
\node[vertex](v1) at (0,2){};
\node[vertex](v2) at (1,1){};
\node[vertex](v3) at (1,0){};
\node[vertex](v7) at (1,-1){};
\node[vertex](v4) at (-1,1){};
\node[vertex](v5) at (0,1){};
\node[vertex](v6) at (0,0){};

\draw[edge](v1)--(v2);
\draw[edge](v2)--(v3);
\draw[edge](v3)--(v7);
\draw[edge](v1)--(v4);
\draw[edge](v1)--(v5);
\draw[edge](v5)--(v6);
\end{scope}

\begin{scope}[shift= {(-2.5,-0.5)}]
\node[vertex](v1) at (0,2){};
\node[vertex](v2) at (-1,1){};
\node[vertex](v4) at (0,1){};
\node[vertex](v5) at (1,1){};

\draw[edge](v1)--(v2);
\draw[edge](v1)--(v4);
\draw[edge](v1)--(v5);
\end{scope}

\end{scope}
\end{tikzpicture}

%% file: f-bonsma-clawfree.tex
\begin{tikzpicture}
\begin{scope}[scale = 0.8]
\begin{scope}

\foreach	\i / \j / \k in {0/0/2,  1/0/3,  2/1/2,  3/1/3}{
	\node[rvertex](a\i) at (\j, \k){};
}
\foreach	 \i / \j in {0/1,  0/2,  0/3,  1/2,  1/3,  2/3}{
	\draw[edge](a\i) -- (a\j);
}

\foreach	\i / \j / \k in {0/3/2,  1/3/3,  2/4/2,  3/4/3}{
	\node[rvertex](b\i) at (\j, \k){};
}
\foreach	 \i / \j in {0/1,  0/2,  0/3,  1/2, 1/3,  2/3}{
	\draw[edge](b\i) -- (b\j);
}

\foreach	\i / \j / \k in {0/1.3/0,  1/1.3/1,  2/2.7/0,  3/2.7/1,  4/2/0.5}{
	\node[rvertex](c\i) at (\j, \k){};
}
\foreach	 \i / \j in {0/1,  0/4,  1/4,  2/3, 2/4,  3/4}{
	\draw[edge](c\i) -- (c\j);
}

\node[rvertex](u1) at (0.5, 1){};
\node[rvertex](u2) at (3.3, 1.3){};
\node[bvertex](u3) at (5/3, 3.2){};
\node[bvertex](u4) at (7/3,3.2){};

\draw[edge](u1)--(a0);
\draw[edge](u1)--(c0);
\draw[tedge](u3)--(a3);
\draw[edge](u3)--(u4);
\draw[tedge](u4)--(b1);
\draw[edge](u2)--(c3);
\draw[edge](u2)--(b2);
\draw[edge](a2)--(b0);
\draw[edge](c2) to [bend  right = 50](b3);

\node[](t) at (2, -1){$G_1$};

\end{scope}

\begin{scope}[shift= {(6,0)}]

\foreach	\i / \j / \k in {0/0/2,  1/0/3,  2/1/2,  3/1/3}{
	\node[bvertex](a\i) at (\j, \k){};
}
\foreach	 \i / \j in {0/1,  0/2,  0/3,  1/2,  1/3,  2/3}{
	\draw[edge](a\i) -- (a\j);
}

\foreach	\i / \j / \k in {0/3/2,  1/3/3,  2/4/2,  3/4/3}{
	\node[bvertex](b\i) at (\j, \k){};
}
\foreach	 \i / \j in {0/1,  0/2,  0/3,  1/2,  1/3,  2/3}{
	\draw[edge](b\i) -- (b\j);
}

\foreach	\i / \j / \k in {0/1.3/0,  1/1.3/1,  2/2.7/0,  3/2.7/1,  4/2/0.5}{
	\node[rvertex](c\i) at (\j, \k){};
}
\foreach	 \i / \j in {0/1,  0/4,  1/4,  2/3, 2/4,  3/4}{
	\draw[edge](c\i) -- (c\j);
}

\node[rvertex](u1) at (0.5, 1){};
\node[bvertex](u2) at (3.3, 1.3){};
\node[bvertex](u3) at (2, 3){};

\draw[tedge](u1)--(a0);
\draw[edge](u1)--(c0);
\draw[edge](u3)--(a3);
\draw[edge](u3)--(b1);
\draw[tedge](u2)--(c3);
\draw[edge](u2)--(b2);
\draw[edge](a2)--(b0);
\draw[tedge](c2) to [bend  right = 50](b3);

\node[](t) at (2, -1){$G_2$};

\end{scope}

\begin{scope}[shift={(11.5,0)}]

	\node[vertex, minimum size = 3mm, fill = none](a) at (0.5,2.5){};
	\node[vertex, minimum size = 3mm, fill = none](b) at (3.5, 2.5){};
	\node[vertex, minimum size = 3mm, fill = none](c) at (2, 0.5){};




\draw[edge](a)[bend angle=10, bend left]to (b);
\draw[edge](a)[bend angle=10, bend right]to (b);
\draw[tedge](a)--(c);
\draw[tedge](b)[bend angle=10, bend left]to (c);
\draw[tedge](b)[bend angle=10, bend right]to (c);

\node[](t) at (2, -1){$G_2^*$};

\end{scope}

\end{scope}

\end{tikzpicture}